\documentclass[a4paper,11pt]{amsart}

\usepackage{graphicx}
\usepackage{mathptmx}
\usepackage{amsmath}
\usepackage{amssymb}
\usepackage{enumitem}
\usepackage{xcolor}

\newmuskip\pFqmuskip

\newcommand*\pFq[6][8]{%
  \begingroup 
  \pFqmuskip=#1mu\relax
  \mathcode`=\string"8000
  \begingroup\lccode`\~=`\,
  \lowercase{\endgroup\let~}\pFqcomma
  F^{#2}_{#3}{\left(\genfrac..{0pt}{}{#4}{#5}\bigg|#6\right)}%
  \endgroup
}
\newcommand{\pFqcomma}{\mskip\pFqmuskip}

\newtheorem{theorem}{Theorem}[section]

\begin{document}

\title[Spivey's type recurrence relation for Lah-Bell polynomials]{Spivey's type recurrence relation for Lah-Bell polynomials}

\author{Taekyun  Kim}
\address{Department of Mathematics, Kwangwoon University, Seoul 139-701, Republic of Korea}
\email{tkkim@kw.ac.kr}
\author{Dae San  Kim}
\address{Department of Mathematics, Sogang University, Seoul 121-742, Republic of Korea}
\email{dskim@sogang.ac.kr}


\begin{abstract}
The aim of this paper is to derive Spivey's type recurrence relations for the Lah-Bell polynomials and the $r$-Lah-Bell polynomials by utilizing operators $X$ and $D$ satisfying the commutation relation $DX-XD=1$. Here $X$ is the `multiplication by $x$' operator and $D$ is the differentiation operator $D=\frac{d}{dx}$. In addition, we obtain Spivey's type recurrence relation for the $\lambda$ analogue of $r$-Lah-Bell polynomials by some other method without using the operators $X$ and $D$.
\end{abstract}

\maketitle

\section{Introduction}
The Stirling number of the second kind ${n \brace k}$ counts the number of partitions of the set $[n]=\{1,2,\dots,n\}$ into $k$ nonempty disjoint subsets (see [18,21,24,26]). The Bell number $\phi_{n}$ enumerates the number of all partitions of the set $[n]$ into nonempty disjoint subsets and hence we have $\phi_{n}=\sum_{k=0}^{n}{n \brace k}$. The Bell polynomial $\phi_{n}(x)$ is the natural polynomial extension of the Bell number $\phi_{n}$ (see [9,10,22,32]), defined by $\phi_{n}(x)=\sum_{k=0}^{n}{n \brace k} x^{k}$ and also given by Dobinski's formula $\phi_{n}(x)=e^{-x}\sum_{k=0}^{\infty} \frac{k^{n}}{k!}x^{k}$.\par
The Lah number $L(n,k)$ is the number of partitions of the set $[n]$ into $k$ nonempty disjoint linear ordered subsets (see [6,14,27,30,34,37]). The Lah-Bell number $LB_{n}$ is the number of all partitions of the set $[n]$ into nonempty disjoint linearly ordered subsets and hence we have $LB_{n}=\sum_{k=0}^{n}L(n,k)$ (see [1,5,6,12-15,27,29]). The Lah-Bell polynomial $LB_{n}(x)$ is the natural polynomial extension of the Lah-Bell number $LB_{n}$, defined by $LB_{n}(x)=\sum_{k=0}^{n}L(n,k)x^{k}$. \par
Let $r$ be a nonnegative integer. The $r$-Lah number $L^{r}(n,k)$ is the number of partitions of the set $[n]$ into $k$ nonempty disjoint linearly ordered subsets in such a way that the number $1,2,3,\dots,r$ are in distinct nonempty disjoint linearly ordered subsets (see [12,30,34]). The $r$-Lah-Bell number $LB_{n}^{(r)}$ is the number of all partitions of the set $[n]$ into nonempty disjoint linearly ordered subsets in such a way that the number $1,2,3,\dots,r$ are in distinct nonempty disjoint linearly ordered subsets and hence we have $LB_{n}^{(r)}=\sum_{k=0}^{n}L^{r}(n,k)$. The $r$-Lah-Bell polynomial $LB_{n}^{(r)}(x)$ is the natural polynomial extension of the $r$-Lah-Bell number $LB_{n}^{(r)}$, defined by $LB_{n}^{(r)}(x)=\sum_{k=0}^{n}L^{r}(n,k)x^{k}$. \par
The Spivey showed by combinatorial method the following recurrence relation for the Bell numbers which is given by
\begin{displaymath}
\phi_{n+m}=\sum_{j=0}^{m}\sum_{k=0}^{n}\binom{n}{k}{m \brace j}j^{n-k}\phi_{k},\quad (n,m\ge 0),\quad (\mathrm{see}\ [35]).
\end{displaymath}
Subsequent to Spivey's work, a lot of researchers have investigated and expanded upon his formula for the Bell numbers (see [20]). \par
Let $X$ be the `multiplication by $x$' operator, and let $D=\frac{d}{dx}$. Then $X$ and $D$ satisfy the commutation relation $DX-XD=1$ (see [3,4,11,16,17,23,25,31]). In this paper, we show Spivey's type recurrence relations for the Lah-Bell polynomials $LB_{n}(x)$ (see [12,13,15,28,29]) and the $r$-Lah-Bell polynomials $LB_{n}^{(r)}(x)$ by using the operators $X$ and $D$. Indeed, for the $r$-Lah-Bell polynomials we show \par
\begin{equation*}
\mathrm{LB}_{n+m}^{(r)}(x)=\sum_{k=0}^{m}\sum_{l=0}^{n}\binom{n}{l}L^{r}(m,k)\langle m+k\rangle_{n-l}x^{k} \mathrm{LB}_{l}^{(r)}(x).
\end{equation*}
Futhermore, we consider the $\lambda$-analogues of $r$-Lah numbers and $r$-Lah-Bell polynomials and show by utilizing some other method without using the operators $X$ and $D$ the following Spivey's type recurrence relation for the $\lambda$-analogue of $r$-Lah-Bell polynomials
\begin{equation*}
\mathrm{LB}_{n+m,\lambda}^{(r)}(x)=\sum_{j=0}^{m}\sum_{k=0}^{n}\binom{n}{k}L_{\lambda}^{r}(m,j)\langle m+j\rangle_{n-k}x^{j}\lambda^{n-k}\mathrm{LB}_{k,\lambda}^{(r)}(x).
\end{equation*} \par
Along the way, we derive the generating function and an explicit expression for the $r$-Lah numbers. We obtain the generating function and Dobinski-like formula for the $r$-Lah-Bell polynomials. We deduce the generating function and an explicit expression for the $\lambda$-analogue of $r$-Lah numbers. We get the generating function and Dobinski-like formula for the $\lambda$-analogue of the $r$-Lah-Bell polynomials. We let the reader refer to [1,2,5,7,8,19,33,36] as general references. For the rest of this section, we recall what are necessary throughout this paper.

\vspace{0.1in}

Let $(x)_{n}$ be the falling factorial sequence with
\begin{displaymath}
(x)_{0}=1,\ (x)_{n}=x(x-1)(x-2)\cdots(x-n+1),\ (n\ge 1),
\end{displaymath}
and let $\langle x\rangle_{n}$ be the rising factorial sequence with
\begin{displaymath}
\langle x\rangle_{0}=1,\quad \langle x\rangle_{n}=x(x+1)(x+2)\cdots(x+n-1),\ (n\ge 1).
\end{displaymath}
For $n,k\ge 0$, the Lah numbers are given by
\begin{equation}
\langle x\rangle_{n}=\sum_{k=0}^{n}L(n,k)(x)_{k},\quad (n\ge 0),\quad (\mathrm{see}\ [12,13]).\label{1}
\end{equation}
From \eqref{1}, we note that
\begin{align}
\Big(\frac{1}{1-t}\Big)^{x}&=\sum_{n=0}^{\infty}\binom{x+n-1}{n}t^{n}=\sum_{n=0}^{\infty}\langle x\rangle_{n}\frac{t^{n}}{n!}=\sum_{n=0}^{\infty}\sum_{k=0}^{n}L(n,k)(x)_{k}\frac{t^{n}}{n!} \label{2}\\
&=\sum_{k=0}^{\infty}\sum_{n=k}^{\infty}L(n,k)\frac{t^{n}}{n!}(x)_{k}. \nonumber
\end{align}
On the other other hand, by binomial expansion, we get
\begin{equation}
\Big(\frac{1}{1-t}\Big)^{x}=\Big(1+\frac{t}{1-t}\Big)^{x}=\sum_{k=0}^{\infty}\frac{1}{k!}\Big(\frac{t}{1-t}\Big)^{k}(x)_{k}. \label{3}	
\end{equation}
By \eqref{2} and \eqref{3}, we get
\begin{equation}
\frac{1}{k!}\Big(\frac{t}{1-t}\Big)^{k}=\sum_{n=k}^{\infty}L(n,k)\frac{t^{n}}{n!},\quad (k\ge 0),\quad (\mathrm{see}\ [12,13]). \label{4}
\end{equation}
From \eqref{4}, we have
\begin{equation}
L(n,k)=\frac{n!}{k!}\binom{n-1}{k-1},\quad (n,k\ge 0),\quad (\mathrm{see}\ [12,13,27,29]).\label{5}
\end{equation} \par
The Lah-Bell polynomials are defined by
\begin{equation}
\mathrm{LB}_{n}(x)=\sum_{k=0}^{n}L(n,k)x^{k},\quad (n\ge 0),\quad (\mathrm{see} \ [13]). \label{6}
\end{equation}
By \eqref{4} and \eqref{6}, we get
\begin{equation}
e^{x(\frac{1}{1-t}-1)}=\sum_{n=0}^{\infty}\mathrm{LB}_{n}(x)\frac{t^{n}}{n!},\quad (\mathrm{see}\ [13]). \label{7}
\end{equation}
From \eqref{7}, we have the Dobinski-like formula for the Lah-Bell polynomials
\begin{equation}
\mathrm{LB}_{n}(x)=e^{-x}\sum_{k=0}^{\infty}\frac{\langle k\rangle_{n}}{k!}x^{k},\quad (n\ge 0). \label{8}
\end{equation} \par
Now, we define the operators $X$ and $D$ by
\begin{equation}
Xf(x)=xf(x),\quad Df(x)=\frac{df(x)}{dx},\quad (\mathrm{see}\ [3,4,9,10,11,17]), \label{9}
\end{equation}
satisfying the commutation relation $DX-XD=1$. \par
From \eqref{9}, we note that
\begin{align}
DX^{k}-X^{k}D&=\big(DX^{k-1}-X^{k-1}D\big)X+X^{k-1}\big(DX-XD\big) \label{10}\\
&=\big(DX-XD\big)X^{k-1}+(k-1)X^{k-1}=kX^{k-1},\quad (k\in\mathbb{N}). \nonumber	
\end{align}
Thus, by \eqref{10}, we get
\begin{equation}
\big(XD\big)X^{k}=X^{k}\big(XD+k\big),\quad (k\in\mathbb{N}),\quad (\mathrm{see}\ [3,4,11,16,17,23,25,31]).\label{11}
\end{equation}

\section{Spivey's type recurrence relations for Lah-Bell and $r$-Lah-Bell polynomials}
For $f(x)=x^{m},\ (m \ge 0)$, we have
\begin{align}
\langle XD\rangle_{n}f(x)&=\big(XD\big)\big(XD+1\big)\big(XD+2\big)\cdots\big(XD+n-1\big)x^{m} \label{12}\\
&=m(m+1)(m+2)\cdots(m+n-1\big)x^{m}\nonumber\\
&=\langle m\rangle_{n} x^{m}=\sum_{k=0}^{n}L(n,k)(m)_{k}x^{m}\nonumber\\
&=\sum_{k=0}^{n}L(n,k)X^{k}D^{k}x^{m}=\sum_{k=0}^{n}L(n,k)X^{k}D^{k}f(x). \nonumber	
\end{align}
Thus, by \eqref{12}, we get
\begin{equation}
\langle XD\rangle_{n}=\sum_{k=0}^{n}L(n,k)X^{k}D^{k},\quad (n\ge 0). \label{13}	
\end{equation}
From \eqref{13}, we note that
\begin{align}
\frac{1}{e^{x}}\langle XD\rangle_{n}e^{x}&=\frac{1}{e^{x}}\sum_{k=0}^{n}L(n,k)X^{k}D^{k}e^{x} \label{14}\\
&=\frac{1}{e^{x}}\sum_{k=0}^{n}L(n,k)x^{k}e^{x}= \mathrm{LB}_{n}(x),\quad (n\ge 0). \nonumber
\end{align}
Therefore, by \eqref{14}, we obtain the following theorem.
\begin{theorem}
For $n\ge 0$, we have
\begin{displaymath}
\frac{1}{e^{x}}\langle XD\rangle_{n}e^{x}= \mathrm{LB}_{n}(x).
\end{displaymath}
\end{theorem}
Now, we observe from \eqref{2} that
\begin{align}
\Big(\frac{1}{1-t}\Big)^{x+y}&=\Big(\frac{1}{1-t}\Big)^{x}\Big(\frac{1}{1-t}\Big)^{y}=\sum_{k=0}^{\infty}\langle x\rangle_{k}\frac{t^{k}}{k!}\sum_{l=0}^{\infty}\langle y\rangle_{l}\frac{t^{l}}{l!} \label{15}\\
&=\sum_{n=0}^{\infty}\sum_{k=0}^{n}\binom{n}{k}\langle x\rangle_{k}\langle y\rangle_{n-l}\frac{t^{n}}{n!},\nonumber	
\end{align}
and that
\begin{align}
\Big(\frac{1}{1-t}\Big)^{x+y}&=\Big(\frac{1}{1-t}\Big)^{x}\Big(\frac{1}{1-t}\Big)^{y} =\sum_{k=0}^{\infty}\binom{x+k-1}{k}t^{k}\sum_{l=0}^{\infty}\binom{y+l-1}{l}t^{l}\label{16}\\
&=\sum_{n=0}^{\infty}\sum_{k=0}^{n}\binom{x+k-1}{k}\binom{y+n-k-1}{n-k}t^{n}. \nonumber
\end{align}
From \eqref{15} and \eqref{16}, we get
\begin{align}
&\langle x+y\rangle_{n}=\sum_{k=0}^{n}\binom{n}{k}\langle x\rangle_{k}\langle y\rangle_{n-k}, \label{17} \\
&\binom{x+y+n-1}{n}=\sum_{k=0}^{n}\binom{x+k-1}{k}\binom{y+n-k-1}{n-k},\quad (n\ge 0). \nonumber
\end{align} \par
For $m,n\ge 0$, we have
\begin{align}
\langle XD\rangle_{n+m}&=XD\big(XD+1\big)\cdots\big(XD+m-1\big)\cdots\big(XD+m+(n-1)\big) \label{18}\\
&=\langle XD\rangle_{m}\langle XD+m\rangle_{n}=\langle XD+m\rangle_{n}\langle XD\rangle_{m}\nonumber\\
&=\langle XD+m\rangle_{n}\sum_{j=0}^{m}L(m,j)X^{j}D^{j}=\sum_{j=0}^{m}L(m,j)\langle XD+m\rangle_{n}X^{j}D^{j}.\nonumber	
\end{align}
From \eqref{11} and \eqref{17}, we note that
\begin{align}
\langle XD+m\rangle_{n}X^{j}&=\big(XD+m\big)\big(XD+m+1\big)\cdots\big(XD+m+n-1)X^{j} \label{19}\\
&=\big(XD+m\big)\cdots\big(XD+m+n-2\big)X^{j}\big(XD+m+n-1+j\big) \nonumber\\
&=\cdots\nonumber\\
&=X^{j}\big(XD+m+j\big)\big(XD+m+j+1\big)\cdots\big(XD+m+j+n-1\big) \nonumber\\
&=X^{j}\langle XD+m+j\rangle_{n}=X^{j}\sum_{k=0}^{n}\binom{n}{k}\langle XD\rangle_{k}\langle m+j\rangle_{n-k}.\nonumber
\end{align}
By \eqref{18} and \eqref{19}, we get
\begin{equation}
\langle XD\rangle_{n+m}=\sum_{j=0}^{m}\sum_{k=0}^{n}L(m,j)\binom{n}{k}\langle m+j\rangle_{n-k}X^{j}\langle XD\rangle_{k}D^{j}. \label{20}	
\end{equation} \par
For $m,n\ge 0$, by \eqref{14} and \eqref{20}, we have
\begin{align}
\mathrm{LB}_{n+m}(x)&=\frac{1}{e^{x}}\langle XD\rangle_{n+m}e^{x}=\frac{1}{e^{x}}\langle XD+m\rangle_{n}\langle XD\rangle_{m}e^{x}\label{21} \\
&=\sum_{j=0}^{m}\sum_{k=0}^{n}\binom{n}{k}L(m,j)\langle m+j\rangle_{n-k}\frac{1}{e^{x}}X^{j}\langle XD\rangle_{k}D^{j}e^{x} \nonumber \\
&=\sum_{j=0}^{m}\sum_{k=0}^{n}\binom{n}{k}L(m,j)\langle m+j\rangle_{n-k}\frac{1}{e^{x}}x^{j}\langle XD\rangle_{k}e^{x} \nonumber \\
&=\sum_{j=0}^{m}\sum_{k=0}^{n}\binom{n}{k}L(m,j)\langle m+j\rangle_{n-k}x^{j} \mathrm{LB}_{k}(x). \nonumber
\end{align}
Therefore, by \eqref{21}, we obtain the following Spivey's type formula for the Lah-Bell polynomials.
\begin{theorem}
For $m,n\ge 0$, we have
\begin{displaymath}
\mathrm{LB}_{n+m}(x)=\sum_{j=0}^{m}\sum_{k=0}^{n}\binom{n}{k}L(m,j)\langle m+j\rangle_{n-k}x^{j} \mathrm{LB}_{k}(x).
\end{displaymath}
In particular, letting $x=1$ gives
\begin{displaymath}
\mathrm{LB}_{n+m}=\sum_{j=0}^{m}\sum_{k=0}^{n}\binom{n}{k}L(m,j)\langle m+j\rangle_{n-k}\mathrm{LB}_{k}.
\end{displaymath}
\end{theorem}
For any integer $r\ge 0$, we consider the $r$-Lah numbers defined by
\begin{equation}
\langle x+r\rangle_{n}=\sum_{k=0}^{n}L^{r}(n,k)(x)_{k},\quad (n\ge 0). \label{22}	
\end{equation}
Thus, by \eqref{2}, \eqref{3} and \eqref{22}, we get
\begin{equation}
\sum_{k=0}^{\infty}\frac{1}{k!}\Big(\frac{1}{1-t}\Big)^{r}\Big(\frac{t}{1-t}\Big)^{k}(x)_{k}=\Big(\frac{1}{1-t}\Big)^{x+r}=\sum_{k=0}^{\infty}\sum_{n=k}^{\infty}L^{r}(n,k)\frac{t^{n}}{n!}(x)_{k}.\label{23}
\end{equation}
From \eqref{23}, we note that
\begin{equation}
\frac{1}{k!}\Big(\frac{1}{1-t}\Big)^{r}\Big(\frac{1}{1-t}-1\Big)^{k}=\sum_{n=k}^{\infty}L^{r}(n,k)\frac{t^{n}}{n!},\quad (k\ge 0). \label{24}	
\end{equation}
By \eqref{24}, we easily get
\begin{equation}
\sum_{n=k}^{\infty}L^{r}(n,k)\frac{t^{n}}{n!}=\frac{t^{k}}{k!}\Big(\frac{1}{1-t}\Big)^{k+r}=\sum_{n=k}^{\infty}\frac{n!}{k!}\binom{n+r-1}{k+r-1}\frac{t^{n}}{n!}.\label{25}
\end{equation}
By comparing the coefficients on both sides of \eqref{25}, we get
\begin{equation}
L^{r}(n,k)=\frac{n!}{k!}\binom{n+r-1}{k+r-1},\quad (n,k,r\ge 0).\label{26}	
\end{equation} \par
In view of \eqref{6}, we define $r$-Lah-Bell polynomials
\begin{equation}
\mathrm{LB}_{n}^{(r)}(x)=\sum_{k=0}^{n}L^{r}(n,k)x^{k},\quad (n\ge 0). \label{27}	
\end{equation}
From \eqref{24} and \eqref{27}, we note that
\begin{equation}
\Big(\frac{1}{1-t}\Big)^{r}e^{x\big(\frac{1}{1-t}-1\big)}=\sum_{n=0}^{\infty} \mathrm{LB}_{n}^{(r)}(x)\frac{t^{n}}{n!}.\label{28}	
\end{equation}
Thus, by \eqref{28}, we get
\begin{align}
\sum_{n=0}^{\infty} \mathrm{LB}_{n}^{(r)}(x)\frac{t^{n}}{n!}&=\frac{1}{e^{x}}e^{\frac{x}{1-t}}\Big(\frac{1}{1-t}\Big)^{r}=\frac{1}{e^{x}}\sum_{k=0}^{\infty}x^{k}\frac{1}{k!}\Big(\frac{1}{1-t}\Big)^{k+r}\label{29}\\
&=\sum_{n=0}^{\infty}\frac{1}{e^{x}}\sum_{k=0}^{\infty}\frac{\langle k+r\rangle_{n}}{k!}x^{k}\frac{t^{n}}{n!}.\nonumber	
\end{align}
By comparing the coefficients on both sides of \eqref{29}, we get the Dobinski-like formula for the $r$-Lah-Bell polynomials
\begin{equation}
\mathrm{LB}_{n}^{(r)}(x)=e^{-x}\sum_{k=0}^{\infty}\frac{\langle k+r\rangle_{n}}{k!}x^{k},\quad (n\ge 0).\label{30}	
\end{equation} \par
For $f(x)=x^{m},\ (m \ge 0)$, we have
\begin{align}
\langle XD+r\rangle_{n}f(x)&=\langle XD+r\rangle_{n}x^{m}\label{31}\\
&=\big(XD+r\big)\big(XD+r+1\big)\cdots\big(XD+r+n-1\big)x^{m}\nonumber\\
&=\langle m+r\rangle_{n}x^{m}=\sum_{k=0}^{n}L^{r}(n,k)(m)_{k}x^{m}\nonumber\\
&=\sum_{k=0}^{n}L^{r}(n,k)X^{k}D^{k}x^{m}=\sum_{k=0}^{n}L^{r}(n,k)X^{k}D^{k}f(x). \nonumber
\end{align}
Thus, by \eqref{31}, we get
\begin{equation}
\langle XD+r\rangle_{n}=\sum_{k=0}^{n}L^{r}(n,k)X^{k}D^{k},\quad (n,r\ge 0). \label{32}	
\end{equation}
From \eqref{32}, we note that
\begin{align}
\frac{1}{e^{x}}\langle XD+r\rangle_{n}e^{x}&=\frac{1}{e^{x}}\sum_{k=0}^{n}L^{r}(n,k)X^{k}D^{k}e^{x}\label{33}\\
&=\frac{1}{e^{x}}\sum_{k=0}^{n}L^{r}(n,k)x^{k}e^{x}\nonumber\\
&=\sum_{k=0}^{n}L^{r}(n,k)x^{k}= \mathrm{LB}_{n}^{(r)}(x),\quad (n\ge 0). \nonumber	
\end{align}
Therefore, by \eqref{33}, we obtain the following theorem.
\begin{theorem}
For $n\ge 0$, we have
\begin{equation}
\frac{1}{e^{x}}\langle XD+r\rangle_{n}e^{x}= \mathrm{LB}_{n}^{(r)}(x). \label{34}
\end{equation}
\end{theorem}
For $m,n\ge 0$, we observe that
\begin{equation}
\langle XD+r\rangle_{n+m}=\langle XD+r\rangle_{m}\langle XD+r+m\rangle_{n}=\langle XD+r+m\rangle_{n}\langle XD+r\rangle_{m}.\label{35}
\end{equation}
For $k\in\mathbb{N}$, by \eqref{11}, we get
\begin{equation}
\langle XD+r+m\rangle_{n}X^{k}=X^{k}\langle XD+r+m+k\rangle_{n}=X^{k}\sum_{l=0}^{n}\binom{n}{l}\langle XD+r\rangle_{l}\langle m+k\rangle_{n-l}. \label{36}
\end{equation}
Therefore, by \eqref{36}, we obtain the following theorem.
\begin{theorem}
For $k\ge 1$, we have
\begin{equation}
\langle XD+r+m\rangle_{n}X^{k}=X^{k}\sum_{l=0}^{n}\binom{n}{l}\langle XD+r\rangle_{l}\langle m+k\rangle_{n-l},\label{37}
\end{equation}
where $m,n$ are nonnegative integers.
\end{theorem}
By \eqref{32}, \eqref{35} and \eqref{37}, we get
\begin{align}
\langle XD+r\rangle_{n+m}&=\langle XD+r+m\rangle_{n}\langle XD+r\rangle_{m} \label{39}\\
&=\langle XD+r+m\rangle_{n}\sum_{k=0}^{m}L^{r}(m,k)X^{k}D^{k}\nonumber\\
&=\sum_{k=0}^{m}L^{r}(m,k)\langle XD+r+m\rangle_{n}X^{k}D^{k}\nonumber\\
&=\sum_{k=0}^{m}\sum_{l=0}^{n}\binom{n}{l}L^{r}(m,k)\langle m+k\rangle_{n-l}X^{k}\langle XD+r\rangle_{l}D^{k}.\nonumber
\end{align}
Therefore, by \eqref{39}, we obtain the following theorem.
\begin{theorem}
For $m,n\ge 0$, we have
\begin{equation}
\langle XD+r\rangle_{n+m}=\sum_{k=0}^{m}\sum_{l=0}^{n}\binom{n}{l}L^{r}(m,k)\langle m+k\rangle_{n-l}X^{k}\langle XD+r\rangle_{l}D^{k}. \label{40}
\end{equation}
\end{theorem}
For $m,n\ge 0$, by \eqref{34} and \eqref{40}, we get
\begin{align}
\mathrm{LB}_{n+m}^{(r)}(x)&=\frac{1}{e^{x}}\langle XD+r\rangle_{n+m}e^{x} \label{41} \\
&=\frac{1}{e^{x}}\sum_{k=0}^{m}\sum_{l=0}^{n}\binom{n}{l}L^{r}(m,k)\langle m+k\rangle_{n-l}X^{k}\langle XD+r\rangle_{l}D^{k}e^{x} \nonumber\\
&=\sum_{k=0}^{m}\sum_{l=0}^{n}\binom{n}{l}L^{r}(n,k)\langle m+k\rangle_{n-l}x^{k}\frac{1}{e^{x}}\langle XD+r\rangle_{l}e^{x}\nonumber\\
&=\sum_{k=0}^{m}\sum_{l=0}^{n}\binom{n}{l}L^{r}(m,k)\langle m+k\rangle_{n-l}x^{k} \mathrm{LB}_{l}^{(r)}(x).\nonumber
\end{align}
Therefore, by \eqref{41}, we obtain the following Spivey's type formula for the $r$-Lah-Bell polynomials.
\begin{theorem}
For $m,n\ge 0$, we have
\begin{equation*}
\mathrm{LB}_{n+m}^{(r)}(x)=\sum_{k=0}^{m}\sum_{l=0}^{n}\binom{n}{l}L^{r}(m,k)\langle m+k\rangle_{n-l}x^{k} \mathrm{LB}_{l}^{(r)}(x).
\end{equation*}
In particular, letting $x=1$ yields
\begin{equation*}
\mathrm{LB}_{n+m}^{(r)}=\sum_{k=0}^{m}\sum_{l=0}^{n}\binom{n}{l}L^{r}(m,k)\langle m+k\rangle_{n-l} \mathrm{LB}_{l}^{(r)}.
\end{equation*}
\end{theorem}

\section{Spivey's type recurrence relation for $\lambda$-analogue of $r$-Lah-Bell polynomials}
Let $\lambda$ be any nonzero real number. Then the degenerate falling factorial sequence $(x)_{n,\lambda}$, and the degenerate rising factorial sequence $\langle x\rangle_{n,\lambda}$ are respectively defined by
\begin{equation}
(x)_{0,\lambda}=1,\quad (x)_{n,\lambda}=x(x-\lambda)(x-2\lambda)\cdots\big(x-(n-1)\lambda\big),\ (n\ge 1),\label{42}
\end{equation}
and
\begin{equation*}
\langle x\rangle_{0,\lambda}=1,\ \langle x\rangle_{n,\lambda}=x(x+\lambda)(x+2\lambda)\cdots\big(x+(n-1)\lambda\big),\ (n\ge 1),\ (\mathrm{see}\ [12,16,17,23]).
\end{equation*}
The degenerate exponentials are given by
\begin{equation}
e_{\lambda}^{x}(t)=\sum_{k=0}^{\infty}(x)_{k,\lambda}\frac{t^{k}}{k!},\quad (\mathrm{see}\ [16,17,23,23,27]). \label{43}
\end{equation} \par
In view of \eqref{22}, we consider the $\lambda$-analogue of $r$-Lah numbers defined by
\begin{equation}
\langle x+r\rangle_{n,\lambda}=\sum_{k=0}^{n}L_{\lambda}^{r}(n,k)(x)_{k,\lambda},\quad (n,r\ge 0).\label{44}
\end{equation}
From \eqref{44}, we note that
\begin{align}
\Big(\frac{1}{1-\lambda t}\Big)^{\frac{x+r}{\lambda}}&=\sum_{n=0}^{\infty}\langle x+r\rangle_{n,\lambda}\frac{t^{n}}{n!}=\sum_{n=0}^{\infty}\sum_{k=0}^{n}L_{\lambda}^{r}(n,k)(x)_{k,\lambda}\frac{t^{n}}{n!}\label{45}\\
&=\sum_{k=0}^{\infty}\sum_{n=k}^{\infty}L_{\lambda}^{r}(n,k)\frac{t^{n}}{n!}(x)_{k,\lambda}.\nonumber
\end{align}
On the other hand, $\big(\frac{1}{1-\lambda t}\big)^{\frac{n+r}{\lambda}}$ is given by
\begin{align}
\Big(\frac{1}{1-\lambda t}\Big)^{\frac{x+r}{\lambda}}&=\Big(\frac{1}{1-\lambda t}\Big)^{\frac{r}{\lambda}}\Big(\frac{1}{1-\lambda t}\Big)^{\frac{x}{\lambda}}=\Big(\frac{1}{1-\lambda t}\Big)^{\frac{r}{\lambda}}\Big(\lambda\frac{t}{1-\lambda t}+1\Big)^{\frac{x}{\lambda}}\label{46}\\
&=\Big(\frac{1}{1-\lambda t}\Big)^{\frac{r}{\lambda}}\sum_{k=0}^{\infty}\frac{1}{k!}\Big(\frac{t}{1-\lambda t}\Big)^{k}(x)_{k,\lambda}. \nonumber
\end{align}
By \eqref{45} and \eqref{46}, we obtain the following theorem.
\begin{theorem}
For $k\ge 0$, we have
\begin{equation}
\frac{1}{k!}\Big(\frac{t}{1-\lambda t}\Big)^{k}\Big(\frac{1}{1-\lambda t}\Big)^{\frac{r}{\lambda}}=\sum_{n=k}^{\infty}L_{\lambda}^{r}(n,k)\frac{t^{n}}{n!}. \label{47}
\end{equation}
\end{theorem}
From \eqref{47}, we have
\begin{align}
\sum_{n=k}^{\infty}L_{\lambda}^{r}(n,k)\frac{t^{n}}{n!}&=\frac{t^{k}}{k!}\Big(\frac{1}{1-\lambda t}\Big)^{\frac{r+\lambda k}{\lambda}}=\frac{t^{k}}{k!}\sum_{n=0}^{\infty}\binom{\frac{r+\lambda k}{\lambda}+n-1}{n}\lambda^{n}t^{n}\label{48}\\
&=\frac{t^{k}}{k!}\sum_{n=0}^{\infty}\frac{\big(r+\lambda(k+n-1)\big)_{n,\lambda}}{n!}t^{n}=\frac{t^{k}}{k!}\sum_{n=0}^{\infty}\binom{r+\lambda(k+n-1)}{n}_{\lambda}t^{n}\nonumber\\
&=\sum_{n=k}^{\infty}\frac{n!}{k!}\binom{r+\lambda(n-1)}{n-k}_{\lambda}\frac{t^{n}}{n!},\nonumber
\end{align}
where $\binom{x}{n}_{\lambda}$ are the degenerate binomial coefficients given by
\begin{equation*}
\binom{x}{n}_{\lambda}=\frac{(x)_{n,\lambda}}{n!},\ (n\ge 1),\quad \binom{x}{0}_{\lambda}=1.
\end{equation*}
Therefore, by \eqref{48}, we obtain the following theorem.
\begin{theorem}
For $n,k\ge 0$, we have
\begin{displaymath}
L^{r}_{\lambda}(n,k)=\frac{n!}{k!}\binom{r+\lambda(n-1)}{n-k}_{\lambda},
\end{displaymath}
where $\binom{x}{k}_{\lambda}$ are the degenerate binomial coefficients.
\end{theorem}
In view of \eqref{27}, we define the $\lambda$-analogue of $r$-Lah-Bell polynomials given by
\begin{equation}
\mathrm{LB}_{n,\lambda}^{(r)}(x)=\sum_{k=0}^{n}L_{\lambda}^{r}(n,k)x^{k},\quad (n\ge 0). \label{49}
\end{equation}
By \eqref{47} and \eqref{49}, we get
\begin{equation}
e^{\frac{x}{\lambda}\big(\frac{1}{1-\lambda t}-1\big)}\Big(\frac{1}{1-\lambda t}\Big)^{\frac{r}{\lambda}}=\sum_{n=0}^{\infty} \mathrm{LB}_{n,\lambda}^{(r)}(x)\frac{t^{n}}{n!}. \label{50}
\end{equation}
Therefore, by \eqref{50}, we obtain the following Dobinski-like theorem.
\begin{theorem}
For $n\ge 0$, we have
\begin{displaymath}
\mathrm{LB}_{n,\lambda}^{(r)}(x)=e^{-\frac{x}{\lambda}}\sum_{k=0}^{\infty}\frac{\langle \lambda k+r\rangle_{n,\lambda}}{k!}\Big(\frac{x}{\lambda}\Big)^{k}.
\end{displaymath}
\end{theorem}
Now, we observe from \eqref{50} that
\begin{align}
&\sum_{m,n=0}^{\infty} \mathrm{LB}_{n+m,\lambda}^{(r)}(t)\frac{x^{n}}{n!}\frac{y^{m}}{m!}=\sum_{m=0}^{\infty}\frac{y^{m}}{m!}D_{x}^{m}\sum_{n=0}^{\infty} \mathrm{LB}_{n,\lambda}^{(r)}(t)\frac{x^{n}}{n!}\label{52}\\
&=\sum_{n=0}^{\infty}\frac{\mathrm{LB}_{n,\lambda}^{(r)}(t)}{n!}\sum_{m=0}^{\infty}\frac{y^{m}}{m!}D_{x}^{m}x^{n}=\sum_{n=0}^{\infty}\frac{\mathrm{LB}_{n,\lambda}^{(r)}(t)}{n!}\sum_{m=0}^{n}\binom{n}{m}y^{m}x^{n-m}\nonumber\\
&=\sum_{n=0}^{\infty}\frac{\mathrm{LB}_{n,\lambda}^{(r)}(t)}{n!}(x+y)^{n}=\bigg(\frac{1}{1-\lambda(x+y)}\bigg)^{\frac{r}{\lambda}}e^{\frac{t}{\lambda}\big(\frac{1}{1-\lambda(x+y)}-1\big)}, \nonumber
\end{align}
where $D_{x}=\frac{d}{dx}$. \par
We note from \eqref{47} and \eqref{50} that
\begin{align}
&\Big(\frac{1}{1-\lambda(x+y)}\Big)^{\frac{r}{\lambda}}e^{\frac{t}{\lambda}\big(\frac{1}{1-\lambda(x+y)}-1\big)}=\bigg(\cfrac{1}{(1-\lambda x)\big(1-\frac{\lambda y}{1-\lambda x}\big)}\bigg)^{\frac{r}{\lambda}}e^{ \frac{t}{\lambda} \big(\frac{1}{(1-\lambda x)(1-\frac{\lambda y}{1-\lambda x})}-1\big)}\label{53}\\
&=\Big(\frac{1}{1-\lambda x}\Big)^{\frac{r}{\lambda}}e^{\frac{t}{\lambda}\big(\frac{1}{1-\lambda x}-1\big)}e^{\frac{t}{\lambda}\big(\frac{1}{1-\lambda x}\big)\big(\frac{1}{1-\frac{\lambda y}{1-\lambda x}}-1\big)} \bigg(\cfrac{1}{1-\frac{\lambda y}{1-\lambda x}}\bigg)^{\frac{r}{\lambda}}	\nonumber\\
&= \Big(\frac{1}{1-\lambda x}\Big)^{\frac{r}{\lambda}}e^{\frac{t}{\lambda}\big(\frac{1}{1-\lambda x}-1\big)}\sum_{j=0}^{\infty}\Big(\frac{t}{\lambda}\Big)^{j}\Big(\frac{1}{1-\lambda x}\Big)^{j}\frac{1}{j!}\bigg(\cfrac{1}{1-\frac{\lambda y}{1-\lambda x}}-1\bigg)^{j}\bigg(\cfrac{1}{1-\frac{\lambda y}{1-\lambda x}}\bigg)^{\frac{r}{\lambda}} \nonumber \\
&= \Big(\frac{1}{1-\lambda x}\Big)^{\frac{r}{\lambda}}e^{\frac{t}{\lambda}\big(\frac{1}{1-\lambda x}-1\big)}\sum_{j=0}^{\infty}t^{j}\sum_{m=j}^{\infty}L_{\lambda}^{r}(m,j)\frac{y^{m}}{m!}\Big(\frac{1}{1-\lambda x}\Big)^{m+j}\nonumber \\
&=\sum_{m=0}^{\infty}\sum_{j=0}^{m}t^{j}L_{\lambda}^{r}(m,j)
\frac{y^{m}}{m!}\sum_{l=0}^{\infty}\langle m+j\rangle_{l}\lambda^{l}\frac{x^{l}}{l!}\sum_{k=0}^{\infty}\mathrm{LB}_{k,\lambda}^{(r)}(t)\frac{x^{k}}{k!}\nonumber \\
&=\sum_{m,n=0}^{\infty}\bigg(\sum_{j=0}^{m}\sum_{k=0}^{n}L_{\lambda}^{r}(m,j)\binom{n}{k}t^{j}\mathrm{LB}_{k,\lambda}^{(r)}(t)\langle m+j\rangle_{n-k}\lambda^{n-k}\bigg)\frac{x^{n}}{n!}\frac{y^{m}}{m!}. \nonumber
\end{align}
Therefore, by \eqref{52} and \eqref{53}, we obtain the following Spivey's type formula for $\lambda$-analogue of $r$-Lah-Bell numbers.
\begin{theorem}
For $m,n\ge 0$, we have
\begin{displaymath}
\mathrm{LB}_{n+m,\lambda}^{(r)}(t)=\sum_{j=0}^{m}\sum_{k=0}^{n}\binom{n}{k}L_{\lambda}^{r}(m,j)\langle m+j\rangle_{n-k}t^{j}\lambda^{n-k}\mathrm{LB}_{k,\lambda}^{(r)}(t).
\end{displaymath}
In particular, letting $t=1$ yields
\begin{displaymath}
\mathrm{LB}_{n+m,\lambda}^{(r)}=\sum_{j=0}^{m}\sum_{k=0}^{n}\binom{n}{k}L_{\lambda}^{r}(m,j)\langle m+j\rangle_{n-k}\lambda^{n-k}\mathrm{LB}_{k,\lambda}^{(r)}.
\end{displaymath}
\end{theorem}

\section{Conclusion}
In this paper, we showed Spivey's type recurrence relations for the Lah-Bell polynomials and the $r$-Lah-Bell polynomials by using the `multiplication by $x$' operator $X$ and the differentiation operator $D=\frac{d}{dx}$, which satisfy the commutation relation $DX-XD=1$. In addition, we obtained Spivey's type recurrence relation for the $\lambda$-analogue of $r$-Lah-Bell polynomials without using the operators $X$ and $D$.\par
Along the way, we derived the generating functions and explicit expressions for the $r$-Lah numbers and the $\lambda$-analogue of $r$-Lah numbers. For the latter, we obtained
\begin{align*}
&\frac{1}{k!}\Big(\frac{t}{1-\lambda t}\Big)^{k}\Big(\frac{1}{1-\lambda t}\Big)^{\frac{r}{\lambda}}=\sum_{n=k}^{\infty}L_{\lambda}^{r}(n,k)\frac{t^{n}}{n!}, \\
&L^{r}_{\lambda}(n,k)=\frac{n!}{k!}\binom{r+\lambda(n-1)}{n-k}_{\lambda}.
\end{align*}
In addition, we deduced the generating functions and Dobinski-like formulas for $r$-Lah-Bell polynomials and $\lambda$-analogue of $r$-Lah-Bell polynomials. For the latter, we got
\begin{align*}
&e^{\frac{x}{\lambda}\big(\frac{1}{1-\lambda t}-1\big)}\Big(\frac{1}{1-\lambda t}\Big)^{\frac{r}{\lambda}}=\sum_{n=0}^{\infty} \mathrm{LB}_{n,\lambda}^{(r)}(x)\frac{t^{n}}{n!}, \\
&\mathrm{LB}_{n,\lambda}^{(r)}(x)=e^{-\frac{x}{\lambda}}\sum_{k=0}^{\infty}\frac{\langle \lambda k+r\rangle_{n,\lambda}}{k!}\Big(\frac{x}{\lambda}\Big)^{k}.
\end{align*}

\end{document}